\documentclass{article}
\usepackage{amsfonts}
\makeatletter
  \def\restrict#1{\ifmmode\def\next{\mathpalette
  \mathrestrict}\fi\left.\next{#1}\right._}
  \def\mathrestrict#1#2{\setbox0=\hbox{$\m@th#1{#2}$}\relax
  \dimen@=\dp0 \advance\dimen@ by 0.7ex\relax
  #2\,\vrule height\ht0 width 0.4pt depth\dimen@}
\makeatother
\makeatletter
\@addtoreset{equation}{section}
\makeatother

\makeatletter
\def\@thmcounterend{.}

\def\newtheorem{\@ifstar{\@sthm}{\@Sthm}}

\def\@Sthm#1{\@ifnextchar[{\@othm{#1}}{\@nthm{#1}}}

\def\@xnthm#1#2[#3]#4#5{\expandafter\@ifdefinable\csname
#1\endcsname
   {\@definecounter{#1}\@addtoreset{#1}{#3}%
   \expandafter\xdef\csname the#1\endcsname{\expandafter\noexpand
     \csname the#3\endcsname \@thmcountersep \@thmcounter{#1}}%
   \expandafter\xdef\csname #1name\endcsname{#2}%
   \global\@namedef{#1}{\@thm{#1}{\csname
#1name\endcsname}{#4}{#5}}%
                             
\global\@namedef{end#1}{\@endtheorem}}}

\def\@ynthm#1#2#3#4{\expandafter\@ifdefinable\csname #1\endcsname
   {\@definecounter{#1}%
   \expandafter\xdef\csname the#1\endcsname{\@thmcounter{#1}}%
   \expandafter\xdef\csname #1name\endcsname{#2}%
   \global\@namedef{#1}{\@thm{#1}{\csname
#1name\endcsname}{#3}{#4}}%
                              
\global\@namedef{end#1}{\@endtheorem}}}

\def\@othm#1[#2]#3#4#5{%
  \@ifundefined{c@#2}{\@latexerr{No theorem environment `#2'
defined}\@eha}%
  {\expandafter\@ifdefinable\csname #1\endcsname
  {\global\@namedef{the#1}{\@nameuse{the#2}}%
  \expandafter\xdef\csname #1name\endcsname{#3}%
  \global\@namedef{#1}{\@thm{#2}{\csname
#1name\endcsname}{#4}{#5}}%
  \global\@namedef{end#1}{\@endtheorem}}}}

\def\@thm#1#2#3#4{\refstepcounter
{#1}\@ifnextchar[{\@ythm{#1}{#2}{#3}{#4}}{\@xthm{#1}{#2}{#3}{#4}}}

\def\@xthm#1#2#3#4{\@begintheorem{#2}{\csname
the#1\endcsname}{#3}{#4}%
                    \ignorespaces}

\def\@ythm#1#2#3#4[#5]{\@opargbegintheorem{#2}{\csname
       the#1\endcsname}{#5}{#3}{#4}\ignorespaces}

\def\@begintheorem#1#2#3#4{\trivlist
                 \item[\hskip\labelsep{#3#1\ #2\@thmcounterend}]#4}

\def\@opargbegintheorem#1#2#3#4#5{\trivlist
      \item[\hskip\labelsep{#4#1\ #2\ (#3)\@thmcounterend}]#5}

\def\@sthm#1#2{\@Ynthm{#1}{#2}}

\def\@Ynthm#1#2#3#4{\expandafter\@ifdefinable\csname #1\endcsname
   {\global\@namedef{#1}{\@Thm{\csname #1name\endcsname}{#3}{#4}}%
    \expandafter\xdef\csname #1name\endcsname{#2}%
    \global\@namedef{end#1}{\@endtheorem}}}

\def\@Thm#1#2#3{\@ifnextchar[{\@Ythm{#1}{#2}{#3}}{\@Xthm{#1}{#2}{
#3}}}

\def\@Xthm#1#2#3{\@Begintheorem{#1}{#2}{#3}\ignorespaces}

\def\@Ythm#1#2#3[#4]{\@Opargbegintheorem{#1}
       {#4}{#2}{#3}\ignorespaces}

\def\@Begintheorem#1#2#3{#3\trivlist
                          
\item[\hskip\labelsep{#2#1\@thmcounterend}]}

\def\@Opargbegintheorem#1#2#3#4{#4\trivlist
      \item[\hskip\labelsep{#3#1\ (#2)\@thmcounterend}]}

{\it}{\rm}
{\it}{\rm}
{\bf}{\it}
{\bf}{\it}
{\it}{\rm}
{\bf}{\rm}
{\bf}{\it}
{\it}{\rm}
{\bf}{\rm}
{\it}{\rm}
\newtheorem{prop}{PROPOSITION}{\bf}{\it}
{\it}{\rm}
{\bf}{\rm}
\newtheorem{theorem}{THEOREM}[section]{\bf}{\it}
\newtheorem*{theorem*}{THEOREM}{\bf}{\it}
\newtheorem*{proof}{Proof}{\bf}{\rm}
\newtheorem{remark}{Remark}{\bf}{\rm}
\makeatother
\newcommand{\ze}{\zeta}

\newcommand{\ds}{\displaystyle}

\newcommand{\de}{\delta}

\newcommand{\al}{\alpha}
\newcommand{\be}{\beta}
\newcommand{\ga}{\gamma}

\newcommand{\om}{\omega}

\newcommand{\si}{\sigma}
\newcommand{\Si}{\Sigma}
\newcommand{\nab}{\nabla}
\newcommand{\half}{\frac{1}{2}}

\newcommand{\p}{\partial}

\newcommand{\diag}{\mbox{diag}}

\newcommand{\beq}{\begin{equation}}
\newcommand{\eeq}{\end{equation}}
\newcommand{\bna}{\begin{eqnarray}}
\newcommand{\ena}{\end{eqnarray}}
\newcommand{\ba}{\begin{eqnarray*}}
\newcommand{\ea}{\end{eqnarray*}}
\newcommand{\ben}{\begin{enumerate}}
\newcommand{\een}{\end{enumerate}}
\newcommand{\bi}{\begin{itemize}}
\newcommand{\ei}{\end{itemize}}

\newcommand{\RR}{{\mathbb R}}

\newcommand{\eqref}[1]{{\rm (\ref{#1})}}
\newcommand{\mb}[1]{\mathbf{#1}}

\begin{document}

\title{Formation of Singularities in\\
 Relativistic Fluid Dynamics and in\\
 Spherically Symmetric Plasma Dynamics}

\author{
YAN GUO\thanks{Yan Guo is an Alfred~P.~Sloan research fellow and was supported in part by a National Science
 Foundation Postdoctoral Fellowship and the National Science Foundation
 grant DMS-9623253.}\\
Division of Applied Mathematics\\ Brown University
\and
A. SHADI TAHVILDAR-ZADEH\thanks{A. Shadi Tahvildar-Zadeh  is an Alfred~P.~Sloan research fellow and was supported in part by the National
Science Foundation grant DMS-9704430.}\\
Department of Mathematics\\ Princeton University
}
\date{January 1999}


\maketitle

\section{Introduction} Quasilinear hyperbolic systems have a special
place in the theory of partial differential equations since most  of
the PDEs arising in continuum physics are of this form.  Well-known
examples are the Euler equations for a perfect compressible fluid, the
equations of elastodynamics for a perfect elastic solid, and equations
describing a variety of field-matter interactions, such as
magnetohydrodynamics etc.  It is well-known that for all these  systems
the Cauchy problem is well-posed, i.e., it has a unique classical
solution in a small neighborhood (in space-time) of the hypersurface on
which the initial data are given.

On the other hand, it is not expected that these systems will have a
global-in-time regular solution, because shock discontinuities are
expected to form at some point, at least as long as the initial data
are not very small.  In more than one space dimension, there are no
general theorems to that effect however, mainly because in higher
dimensions, the method of characteristics, which is a powerful tool in
one dimension for the study of hyperbolic systems, becomes
intractable.

In 1985 T.~C.~Sideris published a remarkable paper  on the formation of
singularities in three-dimensional compressible fluids \cite{S-fstd},
proving that the classical solution to Euler equations has to break
down in finite time.  His proof was based on studying certain averaged
quantities formed out of the solution, showing that they satisfy
differential inequalities whose solutions have finite life-span.  Such
a technique was already employed by Glassey \cite{G-busc} in the case
of a nonlinear Schr\"odinger equation.  The idea is that by using
averaged quantities one is able to avoid local analysis of solutions.
The same technique was subsequently used to prove other formation of
singularity theorems: for a compressible fluid body surrounded by
vacuum in the nonrelativistic \cite{MUK-sssc} and relativistic
\cite{R-ivps} cases, for the spherically symmetric Euler-Poisson
equations in the attractive \cite{MP-ssss} and repulsive \cite{P-ngse}
cases, for magnetohydrodynamics \cite{R-ofsm}, and for elastodynamics
\cite{T-rnes}.  In this paper we present two more such ``Siderian''
blowup theorems: one in relativistic fluid mechanics, and the other in
plasma dynamics.

\section{Relativistic Fluid Dynamics}
Let $(M,g)$ be the Minkowski spacetime, with $(x^\mu)$, $\mu =
0,\dots,3$ the global coordinate system on $M$ in which  $g_{\mu\nu} =
\diag(-1,1,1,1)$.  We will use the standard convention that Greek
indices run from 0 to 3, while Latin ones run from 1 to 3.  Indices are
raised and lowered using the metric tensor $g$, and all up-and-down
repeated indices are summed over the range.  We also denote $t = x^0$
 and $ \mb{x}=(x^1,x^2,x^3) $.  In the following, we adopt the notation
 and terminology of \cite{C-sgrf} and quote from it some of the basic
 facts regarding relativistic dynamics:

The energy tensor for a relativistic
perfect  fluid is 
\begin{equation}\label{def:T}
 T^{\mu\nu} = (\rho+p)u^\mu u^\nu + p g^{\mu\nu}.
\end{equation}
In this formula,
\begin{enumerate}
\item $u = (u^0,\mb{u})$ is the four-velocity field of the fluid, a unit
 future-directed timelike vectorfield on $M$,
so that $g(u,u) = -1$ and hence
\[
u^0 =\sqrt{1+|\mb{u}|^2}.
\]
We note  that here, unlike the nonrelativistic case considered by
Sideris, {\em all} components of the energy tensor are quadratic
in the velocity.
\item $\rho\geq 0$ is the {\em proper energy density} of the
 fluid, the eigenvalue of $T$ corresponding to the eigenvector $u$.
 It is a function of the (nonnegative) thermodynamic variables $n$,
 the {\em number
 density} and $s$,
the {\em entropy per particle}.  The particular dependence of $\rho$ on
these variables is given by the {\em equation of state}
\begin{equation}\label{eq:state}
\rho = \rho(n,s).
\end{equation}
\item $p\geq 0$ is the fluid  pressure, defined by
\begin{equation}\label{def:p}
p = n\frac{\p \rho}{\p n} - \rho.
\end{equation}
Basic assumptions on the equation of state of a perfect fluid are
\begin{equation}\label{ces}
\frac{\p \rho}{\p n}>0,\quad \frac{\p p}{\p n} > 0,\quad
\frac{\p \rho}{\p s} \geq 0
\mbox{ and } =0 \mbox{ iff }s = 0.
\end{equation}
In particular, these  insure that $\eta$, the speed of sound in the
fluid, is always real:
\[
\eta^2 := \left(\frac{dp}{d\rho}\right)_s.
\]
In addition,  the energy tensor (\ref{def:T}) must satisfy the
positivity condition, which implies that we must  have
\begin{equation}\label{p<rho}
p \leq \rho.
\end{equation}
\end{enumerate}
A typical example of an equation of state is that of a polytropic gas.
A perfect fluid is called {\em polytropic} if the equation of state is
of the form
\begin{equation}\label{polyt}
\rho =  n + \frac{A(s)}{\ga-1} n^{\ga},
\end{equation}
where $1 <\ga<2$ and $A$ is a positive  increasing function   of $s$
(The speed of light is equal to one).  This implies that
 $p = An^{\ga}$ and thus the sound speed $\eta(n,s)$ is determined
as follows:
\[
\eta^2 = \left(\frac{dp}{d\rho}\right)_s =
 \frac{\p p/\p n}{\p \rho/\p n} =
\frac{\ga(\ga-1)An^{\ga-1}}{\ga-1+\ga An^{\ga-1}}.
\]
In particular, the sound speed is increasing with density and is bounded
 above by $\sqrt{\ga-1}$.

The equations of motion for a relativistic perfect fluid are:
\begin{equation}\label{eq:T}
\p_\nu T^{\mu\nu} = 0.
\end{equation}
Moreover, $n=n(x)$  satisfies the {\em continuity equation  }
\begin{equation}\label{eq:n}
\p_\nu (n u^\nu) = 0.
\end{equation}
Given an equation of state (\ref{eq:state}), the system of equations
 (\ref{eq:T}-\ref{eq:n}) provides 5 equations for the 5 unknowns
$n(x)$, $s(x)$ and $\mb{u}(x)$.
The component of (\ref{eq:T}) in the direction of $u$ is
\begin{equation}\label{eq:rho}
u^\nu \p_\nu \rho + (\rho +p) \p_\nu u^\nu = 0.
\end{equation}
As long as the solution is $C^1$, this is equivalent to the
{\em adiabatic condition}
\begin{equation}\label{eq:adia}
u^\nu \p_\nu s = 0.
\end{equation}
The component of (\ref{eq:T}) in the direction orthogonal to $u$ is
\begin{equation}\label{eq:u}
(\rho+p)u^\nu\p_\nu u^\mu + h^{\mu\nu}\p_\nu p = 0,
\end{equation}
where
\[
h_{\mu\nu} := g_{\mu\nu} + u_\mu u_\nu
\]
 is the projection tensor onto the orthogonal complement of
 $u(x)$ in $T_xM$.

Thus the system of equations for a relativistic fluid can be written as
 follows:
\begin{equation}\label{eq:rf}
\left\{\begin{array}{rcl} \p_\nu (n u^\nu) & = & 0, \\
(\rho+p)u^\nu\p_\nu u^\mu + h^{\mu\nu}\p_\nu p & = & 0, \\
u^\nu \p_\nu s & = & 0.
\end{array}\right.
\end{equation}
The {\em Cauchy problem} for a relativistic fluid consists of
specifying the values of $n$, $s$ and $\mb{u}$ on a spacelike
hypersurface $\Si_0$ of $M$,
\begin{equation}\label{data}
\restrict{n}{\Si_0} = n_0,\quad \restrict{s}{\Si_0} = s_0,\quad
\restrict{\mb{u} }{\Si_0} = \mb{u}_0,
\end{equation}
and finding a solution $(n,\mb{u},s)$ to
(\ref{eq:rf},\ref{data}) in a neighborhood of $\Si_0$ in $M$.  In
 particular,   let $\Si_0 = \RR^3 \times\{0\}$ be the hyperplane $t=0$
  in $M$ and suppose the initial data (\ref{data}) correspond to a
smooth compactly supported perturbation of a quiet  fluid filling the
space, i.e.,  assume
\begin{equation}\label{cspqf}
\parbox{4in}{$n_0$, $s_0$ and $\mb{u}_0$ are smooth functions on
 $\RR^3$ and
there are positive constants $R_0$, $\bar{n}$ and $\bar{s}$ such that
outside the
ball $B_{R_0}(0)$ we have
$n_0 = \bar{n}$, $s_0 = \bar{s}$, and $\mb{u}_0 = 0$.}
\end{equation}
Let $\bar{\eta} = \eta(\bar{n},\bar{s})$ be the sound speed in the background 
quiet state.  We then have
\begin{prop}\label{prop:dod}
Any $C^1$ solution of {\rm(\ref{eq:rf},\ref{data},\ref{cspqf})}
will
satisfy
\[
n = \bar{n},\quad s = \bar{s}, \quad \mb{u} = 0,
\]
outside the ball $B_{R(t)}(0)$ where $R(t) = R_0 + \bar{\eta}t$.
\end{prop}
\begin{proof}
It is enough to check that the system (\ref{eq:rf}) can be written in
symmetric hyperbolic form:
\begin{equation}\label{shf}
 A_{ij}^\mu(U)\p_\mu U^j = 0\mbox{ where }A_{ij}^\mu = A_{ji}^\mu
\mbox{ and }A_{ij}^0 \mbox{ is positive definite}.
\end{equation}
This can be accomplished for example by using $p$ instead of $n$ as an
unknown.  By (\ref{ces}), we can think of $n$ as a function of $p$ and
and $s$ and thus of $\rho$ as a function of $p$ and $s$. By (\ref{def:p}) it is then
 easy to see
 that (\ref{eq:rf}) is equivalent to the following system for the
unknowns $U = (p,u,s)$:
\begin{equation} \label{eq:pus}
\left\{\begin{array}{rcl}
\displaystyle\frac{1}{(\rho+p)\eta}u^\nu\p_\nu p + \eta \p_\nu u^\nu
 & = & 0 \\
\eta h^{\mu\nu}\p_\nu p + (\rho+p)\eta u^\nu \p_\nu u^\mu  & = & 0\\
u^\nu \p_\nu s & = & 0.
\end{array}\right.
\end{equation}
Let $\bar{U} = (\bar{p},1,0,0,0,\bar{s})$ denote the constant
background solution to (\ref{eq:pus}).
Let $\bar{\ze}:=(\bar{ \rho}+\bar{ p})\bar {\eta}>0$. We then have
 that the differential operator $P = \bar{A}^\mu\p_\mu$ corresponding
to the linearization of (\ref{shf}) at $\bar{U}$ is symmetric
hyperbolic, with
\[
\bar{A}^0 = A^0(\bar{U}) = \diag(\frac{1}{\bar{\ze}},\bar{\ze},
\bar{\ze},\bar{\ze},\bar{\ze},1),\qquad
\bar{A}^i = A^i(\bar{U})  = \left(\begin{array}{cccc} 0 & 0 &
\bar{\eta} \mb{e}_i^T & 0 \\
0  & 0 & 0  & 0 \\
\bar{\eta}\mb{e}_i & 0 & 0 & 0 \\
0  & 0 & 0  & 0
\end{array}\right).
\]
Once we have this, we can use energy estimates, as in  \cite{S-fssn},
 to conclude the desired domain of dependence statement.
\end{proof}

We now prove that  for large enough initial data, the solution to
 (\ref{eq:rf},\ref{data},\ref{cspqf}) cannot remain $C^1$ for all
 $t>0$.   Such a result was announced in \cite{R-ivps}, but the 
unpublished proof contained an error which invalidated the argument \cite{R-pc}.

First of all, a scaling analysis shows that without loss of generality
we can set
 $R_0 = 1$.  Let
\[
B_t = \{x \in M \ |\ x^0 = t, \ |\mb{x}| \leq R(t) = 1 + \bar{\eta}t\}
\]
denote the time $t$ slice of the range of influence of the data, and
let
\begin{equation}\label{def:Q}
Q(t) := \int _{\RR^3} g_{ij} x^i T^{0j} = \int \mb{x} \cdot \mb{u} u^0
(\rho+p)
\end{equation}
be the total radial momentum of the fluid at time $t$. We then have
\bna
Q'(t) &= &\int g_{ij} x^i \p_0T^{0j} = -\int g_{ij} x^i
\p_kT^{kj}\nonumber \\
& = &
\int g_{ij}(T^{ij} - \bar{T}^{ij}) = \int (\rho+p)|\mb{u}|^2 +3(p -
\bar{p}).\label{q'}
\ena

Let
\begin{equation}\label{def:E}
E = \int_{\RR^3}  T^{00} - \bar{T}^{00}= \int (\rho+p)|\mb{u}|^2 +
\rho - \bar{\rho}
\end{equation}
be the total energy of the perturbation.  By (\ref{eq:T}) it is a
conserved quantity, $E = E_0$.
Our goal is to use $E$ to obtain a differential inequality for $Q$ that
would lead to blowup.

We are going to make two assumptions on the equation of state of the
fluid, which are quite natural from a physical point of view.  First we
note that, as mentioned before,
we can use the pressure $p$ as a thermodynamic variable in place of
$n$. The equation of state of the fluid then has the form $\rho =
\rho(p,s)$.  The two assumptions are:
\begin{itemize}
\item[{\bf (A1)}]
$\rho(p,s)$ is a non-increasing function of $s$, for each $p$.
\item[{\bf (A2)}]
$\eta(p,s)$ is a non-decreasing function of $p$, for each $s$.
\end{itemize}
These two assumptions are in particular satisfied for a polytropic
equation of state (\ref{polyt}).
In order to see that, we observe that $n = (p/A(s))^{1/\ga}$, and from
there we get
\[
\rho (p,s) = \frac{1}{\ga-1}p + \frac{1}{A^{1/\ga}(s)} p^{1/\ga}.
\]
It is then clear that {\bf (A1)} holds.  Moreover
\[
\eta^2(p,s) = \frac{\ga(\ga-1)A^{1/\ga}(s) p^{(\ga-1)/\ga}}{\ga-1+\ga
A^{1/\ga}(s)
p^{(\ga -1)/\ga}}
\]
shows that {\bf (A2)} is  satisfied.

We also make the following assumptions on the initial data:
\begin{itemize}
\item[\bf (D1)]
$\bar{\eta} < \frac{1}{3}$.
\item[\bf (D2)]
$E > 0$.
\item[\bf (D3)]
$s_0(\mb{x}) \geq \bar{s}$ for all $\mb{x} \in B_0$.
\end{itemize}
By (\ref{eq:adia}), the entropy $s$ is constant along the
flow lines, and thus {\bf (D3)} implies that $s(x)\geq \bar{s}$ for
$x\in B_t$.  By {\bf(A1)} and {\bf (A2)} we  then have
\ba
\rho - \bar{\rho} & = & \rho(p,s) - \rho(\bar{p},\bar{s}) = \rho(p,s) -
\rho(p,\bar{s}) + \rho(p,\bar{s}) - \rho(\bar{p},\bar{s}) \leq
\rho(p,\bar{s}) - \rho(\bar{p},\bar{s}) \\
& = & \int_{\bar{p}}^p \frac{\p \rho}{\p p}(p',\bar{s}) \, dp' =
\int_{\bar{p}}^p \frac{1}{\eta^2(p',\bar{s})}\, dp' \leq
\frac{1}{\bar{\eta}^2} (p - \bar{p}).
\ea
By (\ref{q'}) and (\ref{def:E}) we then obtain
\[
Q'(t) \geq 3\bar{\eta}^2 E + (1-3\bar{\eta}^2)\int (\rho+p)|\mb{u}|^2,
\]
which implies, by virtue of {\bf(D1)} and {\bf (D2)} that
\[
Q'(t) \geq (1-3\bar{\eta}^2) \int (\rho+p)|\mb{u}|^2 > 0.
\]
In particular $Q(t)>0$ if $Q(0)>0$.

On the other hand, we can always estimate $Q(t)$ from above, using
(\ref{p<rho}):
\ba
Q^2(t) & \leq & \left( \int(\rho+p)|\mb{u}|^2\right) R^2(t)
 \left(\int_{B_t} (\rho+p)(|\mb{u}|^2+1)\right)\\
& \leq & 2 \left( \int(\rho+p)|\mb{u}|^2\right) R^2(t)
\left(\int_{B_t} (\rho+p)|\mb{u}|^2 + \rho  - \bar{\rho} + \bar{\rho}
 \right)\\
& \leq & \frac{2}{1-3\bar{\eta}^2}Q'(t)R^2(t)[E+\frac{4\pi}{3}\bar{\rho}R^3(t)].
\ea

Integrating this differential inequality and changing the integration 
variable to $r=R(t)$, we obtain
\[
\frac{1}{Q(t)} \leq \frac{1}{Q(0)} -
\frac{1-3\bar{\eta}^2}{2\bar{\eta}}
 \int_{1}^{R(t)} \frac{dr}{Er ^2+ \frac{4\pi}{3}\bar{\rho}r^5},
\]
which contradicts the positivity of $Q$ for all time provided the
initial data satisfies the following  final assumption:
\begin{itemize}
\item[\bf (D4)] $Q(0) >
\displaystyle\frac{2\bar{\eta}}{1-3\bar{\eta}^2}\left( \int_{1}^\infty
 \frac{dr}{E r ^2+ \frac{4\pi}{3}\bar{\rho}r^5} \right)^{-1}$.
\end{itemize}
The contradiction implies that there exists a certain $T^*<\infty$
by which time a $C^1$ solution has to have broken down.  In particular,
the domain of dependence may break down at an earlier time, perhaps
because a shock discontinuity forms.  We have thus proved
\begin{theorem}
Suppose that the equation of state of a  fluid satisfies
 {\bf (A1)} and {\bf (A2)}.
  Then the Cauchy problem
{\rm (\ref{eq:rf},\ref{data},\ref{cspqf})} with initial data
satisfying {\bf (D1--D4)}
 cannot have a
global-in-time $C^1$ solution.
\end{theorem}
\begin{remark}
It is easy to obtain a simpler, sufficient condition for blowup: Let
\[
f(y) : = \left( \int_1^\infty \frac{dr}{r^2 (r^3+y)}\right)^{-1}.
\]
By {\bf (D4)} we thus need
\begin{equation}\label{qbound}
Q(0) > \frac{2\bar{\eta}}{1-3\bar{\eta}^2} \frac{4\pi}{3}\bar{\rho}\;
f(\frac{E}{\frac{4\pi}{3}\bar{\rho}}).
\end{equation}
It is easy to see that $f(0) = 4$, $f'(0) = 16/7$ and that $f$
is a concave function of $y$,  so that
$f(y) < \frac{16}{7} y + 4$.  It is therefore enough to have
\begin{equation}\label{QE}
Q(0) > \frac{32\bar{\eta}}{7(1-3\bar{\eta}^2)} ( E +
\frac{7\pi}{3}\bar{\rho}).
\end{equation}
\end{remark}
We note that unlike the nonrelativistic case, the lower bound for
the initial radial momentum in {\bf(D4)} or (\ref{QE})
depends on the initial energy, and thus on the initial velocity.
  Since $Q$ is of the same order of magnitude as $E$, it is worthwhile
to show that there
{\em exist}  data sets satisfying these largeness conditions. 
 In fact, (\ref{QE}) can be satisfied for $\bar{n}$ small enough. 
 All that is needed is  $\p\rho/\p n>0$ at $n=0$.  We illustrate this in the 
following by considering the polytropic case.

Let us consider a fluid with a polytropic equation of state
(\ref{polyt}), and
consider initial data of the following form
\begin{equation}\label{ex:data}
n_0(\mb{x}) = \bar{n} \psi(r),\qquad \mb{u}_0 (\mb{x}) = \frac{\mb{x}}{r}
 \phi(r),\qquad s_0(\mb{x}) = \bar{s} + \phi(r),
\end{equation}
where $r = |\mb{x}|$.  $\phi$ and $\psi$ are smooth, positive  functions on
$[0,\infty)$ such that
\begin{equation}\label{cond:phi}
\phi(r)\equiv 0\mbox{ for }r\geq 1,\qquad \phi(0) = 0,
\end{equation}
and
\begin{equation}\label{cond:psi}
\psi(r)\equiv 1\mbox{ for }r\geq 1,\quad \int_0^1 (\psi(r)-1) r^2 dr =0.
\end{equation}
We then compute
\ba
E & = & \int_{B_0} (\rho_0+p_0)|\mb{u}_0|^2 +\rho_0 -
\bar{\rho}\\
& = & 4\pi\bar{n}\int_0^1 \left\{\psi\phi^2 +
\frac{1}{\ga-1} \bar{n}^{\ga-1} [A(s)\psi^{\ga}(\ga\phi^2+1)-A(\bar{s})]\right\}
r^2dr,
\ea
and thus $E > 0$ by {\bf (D3)} and (\ref{cond:psi}).  
Now,
\[
Q(0) = 4\pi \bar{n}\int_0^1 \phi\sqrt{1+\phi^2}(\psi +
A\frac{\ga}{\ga-1}
\bar{n}^{\ga-1}\psi^{\ga}) r^3 dr.
\]
Dividing (\ref{QE}) by $4\pi \bar{n}$, all we need  is that the following 
inequality be satisfied for $\bar{n}$ small
enough:
\begin{equation}\label{qineq}
\int_0^1 \psi\phi\sqrt{1+\phi^2 } r^3 dr + O(\bar{n}^{\ga-1}) >
\frac{32}{7(1-3\bar{\eta}^2)} \bar{\eta} \left\{ \int_0^1 \psi\phi^2
r^2dr + \frac{7}{12} +
O(\bar{n}^{\ga-1}) \right\}.
\end{equation}
This is clearly true since  $\bar{\eta} \to 0$ as $\bar{n} \to 0$.  We have thus shown
\begin{prop}
Let $\phi$ and $\psi$ be two smooth, positive functions on $[0,\infty)$
satisfying
{\rm(\ref{cond:phi},\ref{cond:psi})}.  Then there exists $\bar{n}>0$
small enough (depending on $\phi$, $\psi$ and $\ga$) such that the
initial data set $(n_0,\mb{u}_0,s_0)$ of the form {\rm(\ref{ex:data})}
satisfy the conditions {\bf(D1--D4)}, and thus lead
to a blowup for {\rm (\ref{eq:rf},\ref{data},\ref{cspqf})}.
\end{prop}

\section{Euler-Maxwell with Constant Background Charge}

A simple two-fluid model to describe  plasma dynamics is the
so called
Euler-Maxwell system, where a compressible electron fluid interacts
with
a constant ion background. Let $n(t,\mb{x})$, $s(t,\mb{x})$ and  $\mb{v}(t,\mb{x})$ be
the average electron density, entropy, and velocity,
 let $\bar{n}$ be the constant ion density,
and let  $\mb{E}(t,\mb{x})$ and $\mb{B}(t,\mb{x})$ be the electric and
magnetic fields.  Let $c =$  speed of light in vacuum, $e =$ the charge
of an electron, and $m = $ the mass of an electron. The  Euler-Maxwell
system (see \cite[pp. 490--491]{J-ce}) then takes the form:
\begin{equation}\label{eq:em}
\left\{\begin{array}{rcl}
\p_t n + \p_i (nv^i) & = & 0\\
\p_t \Pi^i + \p_j T^{ij}& = & \ds\frac{e\bar{n}}{m} E^i\\
\p_t s + v^i \p_i s & = &  0\end{array}
\qquad
\begin{array}{rcl}
\p_t B^i + c (\nab\times \mb{E})^i & = &  0\\
\p_t E^i - c (\nab\times \mb{B})^i & = & - 4\pi e n v^i ,
\end{array}\right.
\end{equation}
together with the  constraint equations
\begin{equation}\label{eq:const}
\p_i E^i = 4\pi e (n-\bar{n}), \qquad
\p_i B^i = 0.
\end{equation}
In the above, $\Pi$ is the momentum vector,
\[
\Pi = n \mb{v} + \frac{1}{4\pi m c} (\mb{E}\times \mb{B}),
\]
and $T$ is the stress tensor, which can be decomposed into material and
electromagnetic parts:
$T = T_M + T_E$, with
\ba
T_M^{ij} & = & n v^iv^j+\frac{1}{m}p\de^{ij},\\
T_E^{ij} & = & \frac{1}{4\pi m}[\half(|\mb{E}|^2+|\mb{B}|^2)\de^{ij} -
E^iE^j-B^iB^j].
\ea
$p$ is the electron pressure, which is modeled by a polytropic law
$p(n,s)=A(s)n^\ga$,  where $\ga>1$ and $A$ is a positive increasing function.

The system (\ref{eq:em}) being hyperbolic, we once again have the
domain of dependence property.  However, this time the largest
characteristic speed in the background will be $c$, the speed of
light.  We recall that in Sideris's original argument \cite{S-fstd},
the largeness condition on the initial data implied that the initial
velocity had to be supersonic at some point, relative to the sound
speed in the background.  An analogous result in the Euler-Maxwell case
would thus require that the initial velocity be superluminar at some
point, which is absurd.  However, we note that if the data is
spherically symmetric, so will be the solution, and thus there will be
no electromagnetic waves, and the largest characteristic speed will
once again be the sound speed, so a Siderian blowup theorem is possible
in the spherically symmetric case.  Moreover, since spherical symmetry
implies that the flow is irrotational, such a blowup result is
complementary to the recent construction \cite{G-sifl} of global smooth
irrotational solutions with small amplitude for the above system.   We
note that a blowup result in the spherically symmetric, isentropic case
with no background charge has  been obtained \cite{E-fsee}
using Riemann invariants.

\begin{remark}
Under the assumption of spherical symmetry, the Euler-Maxwell system
reduces to what is often referred to as the spherically symmetric
Euler-Poisson system (with repulsive force).  We note the important
distinction  between this, and the general Euler-Poisson system obtained by
taking the Newtonian limit $c \to \infty$ in (\ref{eq:em}).  The latter
is not a hyperbolic system, and does not have finite
propagation speeds.
\end{remark}

We have the following theorem:
\begin{theorem}
Let $\nu_0$, $\si_0$ and  $u_0$ be smooth functions on $\RR^+$ satisfying
\[
u_0(r) \equiv \si_0(r) \equiv \nu_0(r)\equiv 0\mbox{ for }r\geq 1,\qquad u_0(0) = 0,\qquad \si_0(r)\geq 0,
\]
and the {\em neutrality condition}
\begin{equation}
\int_0^1 \nu_0(r) r^2 dr =0.\label{neutral}
\end{equation}
Let $\bar{s} \geq 0$ be fixed.  Then
\begin{itemize}
\item[\rm(a)] There exists $T>0$ and functions $\nu,\si,u,E \in C^1(
[0,T)\times\RR^+)$
such that
\[
\begin{array}{l}\nu(0,r)=\nu_0(r)\\
 \si(0,r) = \si_0(r)\\
  u(0,r)=u_0(r)\end{array}\qquad E(0,r)=\frac{4\pi e}{r^2}
\int_0^r \nu_0(r'){r'}^2 dr'.
\]
and  such that
the Euler-Maxwell
system {\rm(\ref{eq:em})} has a unique solution of the form:
\begin{equation}\begin{array}{rcl}
n(t,\mb{x})&=&\bar{n}+\nu(t,r)  \\
s(t,\mb{x}) &=& \bar{s} + \si(t,r)  \\
\mb{v}(t,\mb{x})&=&u(t,r)\ds\frac{\mb{x}}{r}
\end{array}\
\qquad
\begin{array}{rcl}
\mb{E}(t,\mb{x})& = & E(t,r)\ds\frac{\mb{x}}{r}\\
\mb{B}(t,\mb{x})&\equiv & 0,
\end{array}
\label{sphere}
\end{equation}
where $r = |\mb{x}|$.

\item[\rm(b)] For $t\in [0,T)$, $(n, s, \mb{v} , \mb{E})$ satisfy the
reduced Euler-Maxwell system:
\begin{equation}\label{eq:rem}
\left\{\begin{array}{rcl}
\p_t n + \p_i (nv^i)&=& 0\\
\p_t s + v^i \p_i s &=& 0\\
\p_t (nv^i) + \p_j T^{ij}&=&\ds\frac{e\bar{n}}{m} E^i\\
\p_t E^i  + 4\pi e n v^i & =& 0,
\end{array}
\right.
\end{equation}
where
\[ 
T^{ij} = nv^iv^j +\frac{1}{m} p \de^{ij} + \frac{1}{4\pi m } (\half |\mb{E}|^2 \de^{ij} -  E^i
E^j),
\]
together with the constraint Poisson equation:
\[
\p_i E^i =  4\pi e (n-\bar{n}).
\]

\item[\rm(c)] Let $\bar{\eta} = \sqrt{\ga A(\bar{s})\bar{n}^{\ga-1}}$ be the
sound speed in the background, $R(t) := 1 + \bar{\eta} t$, and let
\[
D_T := \{(t,\mb{x})\ |\ 0\leq t<T,\ |\mb{x}| \geq  R(t) \}.
\]
Then  we have $(n,s,\mb{v}, \mb{E})\equiv (\bar{n},\bar{s},0,0)$ on $D_T$.

\item[\rm(d)] For any fixed $\nu_0(r)$ which satisfies
{\rm(\ref{neutral})}, there exists $u_0(r)$ sufficiently large, such
that the
life-span of the $C^1$ solution {\rm(\ref{sphere})} is finite.
\end{itemize}
\end{theorem}

\begin{proof}
\begin{itemize}
\item[(a)]  The Euler-Maxwell system (\ref{eq:em}) can be written
as a positive, symmetric
hyperbolic system, and therefore has a unique, local $C^1$ solution
with $n>0$
provided  its initial data are sufficiently smooth. Notice that the
initial data are spherically symmetric. Because of the
rotational covariance properties of the Euler-Maxwell system and the
uniqueness of the local solution, the solution remains spherically
symmetric and (a) follows.

\item[(b)] follows since $B^i\equiv 0$.

\item[(c)] Notice that  from the Poisson equation at $t=0$,
\[
E(0,r)=\frac{4\pi e}{r^2}\int_0^r \nu_0(r)r^2 dr \equiv 0\mbox{ for
} r\geq 1
\]
by the
neutrality assumption (\ref{neutral}). Now the reduced
Euler-Maxwell system (\ref{eq:rem}) is still a hyperbolic system, and
we can deduce (c) via the Proposition
 in \cite{S-fssn}.

\item[(d)] Let
\[
Q(t) := \frac{1}{4\pi}\int_{\RR^3} x\cdot \Pi =\int_0^\infty r n u
\,r^2dr.
\]
A direct computation yields:
\[
Q'(t) = \int_0^\infty\left\{nu^2 + \frac{3}{m}(p-\bar{p}) +
\frac{1}{8\pi m}E^2\;
\right\} r^2dr +
\frac{e\bar{n}}{m}\int_0^\infty rE\;r^2dr,
\]
where $\bar{p} = p(\bar{n},\bar{s})$.  Meanwhile, by the first and fourth 
equations in (\ref{eq:rem}),
\[
\int_0^\infty rE(t,r)r^2dr = \int_0^\infty rE(0,r)r^2dr -4\pi e\int_0^t
Q(t')dt'.
\]
Integrating by parts, we notice that
\[
\int_0^\infty rE(0,r)r^2dr=-4\pi e \int_0^\infty  \nu_0(r)r^4dr.
\]
 We now define $\ds y(t) := \int_0^t Q(t')dt'$ and obtain
\begin{equation}\label{ODE}
y''(t) + \om^2 y(t) =  G(t),
\end{equation}
where $\ds\om^2 = \frac{4\pi e^2 \bar{n}}{m}$ is the {\em plasma frequency}, and
\[
G(t) := -\om^2\int_0^\infty  \nu_0(r)r^4dr+ \int_0^\infty
\left\{nu^2 + \frac{3}{m}(p-\bar{p})
 + \frac{1}{8\pi m}E^2\;\right\} r^2dr.
\]
Therefore,
form solving the ODE (\ref{ODE}) for $y(t)$, we have
\begin{equation}
y''(t) = -\om y'(0) \sin \om t + G(t) - \om \int_0^t \sin \om(t-\tau)
G(\tau)d\tau.\label{ode}
\end{equation}
We recall the conserved quantities energy:
\[
\mathcal{E} = \int_0^\infty \left\{\half nu^2+\frac{1}{m(\gamma
-1)}(A(s)n^\gamma -
A(\bar{s})\bar{n}^\gamma)+\frac{1}{8\pi m}E^2\right\} r^2 dr,
\]
and mass
\[
M = \frac{1}{4\pi}\int_{\RR^3} (n - \bar{n}) = \int_0^\infty
\nu(t,r)r^2 dr.
\]
From the neutrality condition
(\ref{neutral}) we have $M \equiv 0$.  Also, $s(0,\mb{x})\geq \bar{s}$ since 
$\si_0\geq 0$ by assumption.  By the adiabatic condition (the second equation 
in (\ref{eq:rem})) entropy is constant along flow lines, and thus $s(t,\mb{x})\geq \bar{s}$ 
for $t<T$.  Since $A$ is an increasing function,
\[
\int A(s) n^\ga - A(\bar{s})\bar{n}^\ga \geq A(\bar{s})\int n^\gamma -\bar{n}^\gamma\ge \bar{\eta}^2 \int n-\bar{n} =0.
\]
Hence we have
\[
\al \mathcal{E} \leq G(t)+
\om^2 \int_0^\infty \nu_0(r)r^4dr\leq \be \mathcal{E},
\]
with $\al = \min\{ 1,3(\ga - 1)\}$, $\be = \max\{2,3(\ga-1)\}$. But for
large
enough $u_0(r)$, $\int_0^\infty \nu_0(r)r^4dr$ is dominated by $\mathcal{E}(0)$.
 Hence, we have
\[
\frac{\al}{2} \mathcal{E} \leq G(t)\leq 2\be \mathcal{E}
\]
for sufficiently large $u_0(r)$. Moreover, we have
\begin{equation}
Q^2(t) \leq R^2(t)\int_0^\infty nu^2 \int_0^{R(t)} n\leq C R^5(t)
\bar{n} \mathcal{E}.
\label{qe}
\end{equation}
 $C$ will henceforth denote a generic numerical constant.  By choosing
 $u_0(r)$ large such that $\mathcal{E}(t)=\mathcal{E}(0)\ge 1$, we have
\[
|y'(0)|=|Q(0)|\le  C \sqrt{\bar{n}} \mathcal{E}.
\]
Thus from
(\ref{ode}), there
 exists
$T_0 = T_0(\ga,\bar{n},\om)>0$ such that for $0\leq t\leq T_0$,
\[
Q'(t) \geq C\mathcal{E}.
\]
Together with (\ref{qe}), we deduce for $0\le t\le T_0$,
$$Q'(t) \geq \frac{C}{R^5(t)\bar{n}}Q^2(t).$$Integrating over $[0,
T_0]$
 we obtain
\begin{equation}
\frac{1}{Q(0)}-\frac{1}{Q(T_0)} \geq \frac{C}{\bar{n}\bar{\eta}}[1-
\frac{1}{(1+\bar{\eta} T_0)^{4}}].\label{qle}
\end{equation}
We can then choose $u_0(r)$ sufficiently large, so that
$Q(0)$ is so large to contradict (\ref{qle}).
\end{itemize}

\end{proof}

\centerline{\bf Acknowledgments}

 We would like to thank S.~Cordier and the anonymous referee for their helpful comments, and  D.~Christodoulou for his  lucid presentation of the results of Sideris, which  provided us  with a  starting point for this project.  

\bibliographystyle{amsplain}
 \ifx\undefined\bysame
\newcommand{\bysame}{\leavevmode\hbox to3em{\hrulefill}\,}
\fi

\end{document}